\documentclass[amsart]{article}

\usepackage{amsmath}
\usepackage{amssymb}
\usepackage{amscd}
\usepackage{bbm}
\usepackage{cyrillic}
\usepackage{latexsym}
\usepackage{makeidx}
\usepackage{mathrsfs}
\usepackage{mathtools}
\usepackage{theorem}
\usepackage[all]{xy}
\usepackage[colorlinks]{hyperref} 
\hypersetup{colorlinks=false , linkcolor=blue}

\DeclareMathAlphabet{\mathbbold}{U}{bbold}{m}{n}

\newtheorem{teo}{Theorem}[section]

\newtheorem{lema}[teo]{Lemma}
\newtheorem{coro}[teo]{Corollary}

{\theorembodyfont{\rmfamily }
\newtheorem{defi}[teo]{Definition}

\newtheorem{ejem}[teo]{Example}
\newtheorem{parra}[teo]{}}
\def\demo{\noindent \textit{Proof: }}
\def\fin{\end{document}}

 \DeclareMathOperator{\hofib}{hofib}
\DeclareMathOperator{\Hom}{Hom} 
 
 \DeclareMathOperator{\rk}{rank}
 
\DeclareMathOperator{\Td}{Td} \DeclareMathOperator{\Th}{Th}

\def\11{1}
\def\11{\mathbbm{1}}

\def\cf{\emph{cf. }}
\def\ch{\mathrm{ch}}

\def\p{\mathfrak{p}}
\def\q{\mathfrak{q}}

\def\raya{\ \underline{\phantom{a}}\ }

\def\t{\mathbf{t}}

\def\qed{\hspace*{\fill }$\square $ }

                    \def\AA{\mathbb{A}}

 \def\EE{\mathbb{E}}

                                         \def\HH{\mathbf{H}}
\def\I{\mathcal{I}}

 \def\NN{\mathbb{N}}
\def\M{\mathcal{M}} \def\MM{\mathbb{M}}
\def\O{\mathcal{O}}
 \def\PP{\mathbb{P}}
                    \def\QQ{\mathbb{Q}}

\def\Z{\mathcal{Z}} \def\ZZ{\mathbb{Z}}

\def\KGL{\mathrm{KGL}}

\def\Gr{\mathrm{Gr}}

\def\SH{\mathbf{SH}}

\def\Sm{\mathbf{Sm}}

\begin{document}

\title{On the higher Riemann-Roch without denominators}

\author{A. Navarro \thanks{Institut f\"{u}r Mathematik, Universit\"{a}t Z\"{u}rich, Switzerland}}

\date{January 25$^{\mathrm{th}}$, 2019}

\maketitle

\begin{abstract}
We prove two refinements of the higher Riemann-Roch without denominators: a statement for regular closed immersions between arbitrary finite dimensional noetherian schemes, with no smoothness assumptions, and a statement for the relative cohomology of a proper morphism.
\end{abstract}

\section*{Introduction}

Grothendieck's Riemann-Roch theorem states that for  $\,f\colon Y\to X\,$  a proper morphism between nonsingular quasi-projective varieties over a field $k$ the square 
$$
\xymatrix{K_0(Y)\ar[r]^{f_!} \ar[d]_{\Td(T_Y)\ch}& K_0(X)\ar[d]^{\Td(T_X)\ch}\\
CH^\bullet (Y) _\QQ \ar[r]^{ f_*} & CH^\bullet (X) _\QQ }
$$
commutes (\emph{cf.} \cite{BorelSerre}). In other words, that for any element $a \in K_0 (Y)$, the formula 
\begin{equation}\label{RRformula}
\Td (T_{\scriptscriptstyle X})\, \ch (f_!(a)) \, =\,  f_* \bigl( \Td
(T_{\scriptscriptstyle Y})\, \ch (a)\bigr) \
\end{equation}
holds in $CH^\bullet (X) _\QQ$. 

In a letter to Serre Grothendieck suggested that in the case  $f$ is a closed immersion the Riemann-Roch should hold in $CH^\bullet (X)$, without the need of tensoring with $\QQ$ (\cf  \cite[0 p. 70]{SGA6}). After the work on the Riemann-Roch at \cite{SGA6} the question still remained as a conjecture known as \emph{the Riemann-Roch without denominators} (\cf \cite[XIV.3]{SGA6}). Jouanolou proved this conjecture in \cite{Jouanolou} for regular closed immersions between quasi-compact schemes over a field, and later for algebraic schemes (\cf \cite[\S 18]{Fulton}). More concretely,  for any $a\in K_0(Y)$ denote $c_i(a)$ its $i$-th Chern class with values in the Chow ring and let $i\colon Y \to X$ be a regular closed immersion of codimension $d$ with normal vector bundle $N$. Jouanolou constructed polynomials with integral coefficients $P_q^d(X;Y)$, such that for any $a\in K_0(Y)$ the formula
$$
c_q(i_*(a))=P_p^d( \rk(a), c_1(a), \ldots , c_{q-d}(a) ; c_1(N), \ldots , c_{q-d}(N))
$$ 
holds in $CH^{q}(X)$. After the coming of higher $K$-theory, Gillet proved in \cite[3]{Gillet} that this same formula holds for elements of higher $K$-theory and $i$ being a closed immersion between smooth schemes over a regular base. Finally, in \cite{KY} transferred Gillet's proof and results into the motivic homotopy context.

In this note we show how to deduce the higher Riemann-Roch without denominators for regular immersions between arbitrary (possibly singular) finite dimensional noetherian schemes from the already known case of regular schemes. More concretely, the theorem we prove in \ref{Teo RR sin deno} is the following:

\medskip

\noindent\textbf{Theorem:} \emph{Let $i\colon Z\to X$ be a regular immersion of codimension $d$. Denote $\bar i_*\colon K(Z)\to K_Z (X)$ as well as  $\bar i_*\colon H_{\mathcal{M}}(Z,\ZZ)\to H_{\mathcal{M},Z}(X,\ZZ) $ the refined Gysin morphisms, $c_{q,r}^Z\colon K_{r,Z}(X)\to H_{\M, Z}(X,\ZZ)$ the $r$-th Chern class with support on $Z$ (cf. \ref{Parra clases Chern Gysin KH y K}  and \ref{Defi Gysin}) and $P_q^d(r,x_1,\cdots ;y_1, \cdots )$ the polynomials with integers coefficients defined in \cite{Jouanolou}. Then for any $a \in KH_r(Z)$ we have}
$$
c_{q,r}^Z(\bar i_*(a))=\bar i_*\bigl(P_q^d(\mathrm{rk}(a), c_{1,r}(a), \ldots, c_{q-d,r}(a); c_1(N_{Z/X}), \ldots ,  c_{q-d} (N_{Z/X}))\bigr).
$$

\medskip

Our proof lies on top of that of Kondo and Yasuda in \cite{KY}. More concretely, all classical computations are summarized in the motivic homotopy context by Kondo and Yasuda in Lemmas \ref{Lema KY 1} and \ref{Lema KY 2}. Our proof relies on a explicit computation of the refined Gysin morphism (that is to say, direct image with support) of the zero section of a vector bundle in \ref{Teo clase de Thom}. Then, the characterization of direct images of \cite{Navarro advances} and the deformation to the normal bundle allows to conclude.

Let $f\colon Y\to X$ be a morphism of schemes. Its inverse image in any reasonable cohomology theory fits into a long exact sequence relating groups called the \emph{relative cohomology} of $f$, which we denote $K(f)$ in the case of $K$-theory. These groups generalize many classical cohomological constructions: the cohomology with proper support, the cohomology with support on a closed subscheme, and the reduced cohomology are the relative cohomology of a closed immersion, an open immersion and the projection over a base point respectively. We prove in Corollary \ref{Coro RR sin deno relativo} a formula without denominators for classes in $K(f)$ for $f$ a proper morphism.

The note is organized as follows: first we recall some basic notation of motivic homotopy theory, in section 2 we describe the Gysin morphism in our context and compute the case of the zero section of a vector bundle, finally we deduce the refinements of the Riemann-Roch.

\bigskip

\noindent\textbf{Acknowledgments:} I would like to thank F. D\'{e}glise, J.I. Burgos Gil, J. Navarro and J. Riou for many helpful comments on a previous version of this text.

\section{Preliminaries}

In this section we recall the notation needed to describe the Gysin morphism constructed in \cite[$\S$ 2]{Navarro advances}. All schemes considered are finite dimensional and noetherian.

\bigskip

Let $X$ be a scheme. Denote $\SH(X)$ the stable homotopy category of Voevodsky (\cite{Voevodsky}). The category $\SH(X)$, whose objects are called \emph{spectra}, is triangulated and monoidal, we denote the unit by $\11_X$ and the $p$-shift functor by  $[p]$ and $\11_X$ respectively. There is also Tate object, $\11_X(1)$ and we denote the $q$-twist functor given  by $(q)$. Hence the $p$-shifted and $q$-twisted spectrum of a $\EE$ of $\SH(X)$ is denoted  $\EE(q)[p]$.

The family of categories $\SH(\raya)$ satisfies Grothendieck's six functor formalism (\cite{Ayoub}). When we consider only schemes over a fixed base $S$, an object $\EE_S$ of $\SH(S)$ (or simply $\EE$ if no confusion is possible) is called an \emph{absolute spectrum} and for $f\colon Y\to S$ we denote $\EE_Y\coloneqq f^*\EE_S$. In general, for any scheme $X$ a spectrum $\EE_X\in \SH(X)$  is also an absolute spectrum \emph{over $X$}.

We say that a spectrum is a \emph{ring} spectrum if it is an associative commutative unitary monoid in $\SH$. Ring spectra define cohomology theories with usual properties. More concretely, for any absolute ring spectrum $\EE$ in $\SH(S)$ we define the $\EE$-cohomology of an $S$-scheme $X$ to be 
$$
\EE^{p,q}(X)\coloneqq \Hom_{\SH(X)}(\11_X, \EE_X(q)[p]).
$$
We also denote $\EE(X)=\bigoplus_{p,q}\EE^{p,q}(X)$.

Let $\EE$ be an absolute ring spectrum, an absolute \emph{$\EE$-module} is a spectrum $\MM$ in $\SH(S)$ together with a morphism of spectra $\upsilon\colon \EE\wedge \MM\to \MM$ in $\SH(S)$ satisfying the usual module condition (\cf \cite[1.1]{Navarro advances}).

For any absolute ring spectrum $\EE$ with unit $e\colon S\to \EE_S $ there is a canonical class in $\eta\in \EE^{2,1}(\PP^1)$ given by the composition $\PP^1\to \11_S(1)[2]\xrightarrow{e\wedge \mathrm{Id}}\EE_S(1)[2]$. We define an \emph{orientation} on $\EE$ to be a class $c_{1} \in \EE ^{2,1} (\PP ^\infty)$ such that for $i_1\colon \PP^1 \hookrightarrow \PP^\infty$ satisfies $i_1 ^* (c_{1})= \eta$. We also say that $\EE$ is \emph{oriented}. Let $(\EE,c_1)$ be an oriented absolute ring spectrum, then for any rank $n$ vector bundle $V\to X$ there are Chern classes $c_i(V)\in \EE^{2i,i}(X)$ for $0<i\leq n$ satisfying the usual Whitney sum and functoriality property. We say that the orientation is \emph{additive} if for two line bundles $L$, $L'$ we have $c_1(L\otimes L')=c_1(L)+c_1(L')$.

\medskip

We denote $\HH(X)$ the unstable homotopy category of Morel and Voevodsky (\cf \cite{MV}), whose objects are called \emph{spaces}. Its pointed version admits a pair of adjoint functors
$$
\Sigma^\infty \colon \HH_\bullet (X)\leftrightarrows \SH(X)\colon \Omega^\infty
$$
By analogy with classical topology, we denote the $\mathrm{Hom}$ sets of $\HH(X)$ by $[\raya, \raya]$.

\begin{ejem}
\begin{enumerate}
\item
The $K$-theory oriented absolute ring spectrum $\KGL$ is defined in \cite{Voevodsky} (\cf also \cite{Riou3}).  It represents Weibel's homotopy invariant $K$-theory for every scheme (\cf \cite[2.15]{Cisinski}), and therefore it represents Quillen's algebraic $K$-theory for regular schemes. We denote the cohomology groups they define as $KH _i (\raya )$. Denote $\Gr$ the infinite Grassmanian, the space $\ZZ\times \Gr$  satisfies that $\Omega^\infty (\KGL)= \ZZ\times \Gr$ and therefore it also represents homotopy invariant $K$-theory. More concretely, denote $S^i$ the space given by the $i$-th simplicial sphere:
$$
KH_i(X)=[S^i\wedge X, \ZZ\times \Gr].
$$

\item
In \cite{Spitz} Spitzweck defined the oriented absolute motivic cohomology ring spectrum $\mathrm{H}_{\Lambda}$ with coefficients in $\Lambda$ for schemes over a Dedekind domain $S$. Motivic cohomology is the universal cohomology with additive orientation, hence there is a cycle class mapping towards more classical cohomologies (such as Betti or algebraic deRham) or towards more refined cohomologies  such as absolute Hodge cohomology with integral coefficients. 

\item
Let $f\colon Y \to X $ be a morphism of schemes and $\EE$ be an absolute ring spectrum. The map $\EE_X\to f_*f^*\EE_Y$ in $\SH(X)$ fits into a distinguished square. Hence, the inverse image in the $\EE$-cohomology fits into a long exact sequence 
$$
\cdots \to \EE^{p-1,q}(Y)\to \EE^{p,q}(f)\to \EE^{p,q}(X)\xrightarrow{f^*} \EE^{p,q}(Y)\to \cdots
$$
where the groups $\EE^{p,q}(f)$ are called \emph{relative $\EE$-cohomology} of $f$ (\cf \cite[\S 1.2]{Navarro advances}). The spectrum representing them is an absolute $\EE$-module over $X$, so we can apply the upcoming Corollary \ref{Coro propiedades Gysin modulos}. However, the spectrum representing relative cohomology has, in general, bad functorial properties. In order to avoid this problem we require $f$ to be proper. We denote by $\hofib(f)$ the space representing these groups in $\HH(S)$, that is the homotopy fiber of $\Omega^\infty\EE_X\to  f_*f^*\Omega^\infty\EE_{X}$. 
\end{enumerate}
\end{ejem}

\section{Gysin morphism}

For the rest of this note, let $(\EE,c_i)$ be an absolute oriented spectrum. 

\begin{parra}\label{Clases superiores con soporte}
The $i$-th Chern class defines a natural transformation
$$
c_i\colon  K_0(\raya ) \to \EE ^{2i,i}(\raya )
$$
of presheaves of sets on $\Sm /S$ which maps every locally free module to its $i$-th Chern class. Riou proved that we have an isomorphism
$$
\Hom(K _0(\raya ), \EE^{2i,i}(\raya ))\simeq [\ZZ\times
\Gr, \Omega^\infty \EE(i)[2i]]
$$
(cf. \cite[1.1.6]{Riou3}). Since $\ZZ\times \Gr$ represents homotopy invariant $K$-theory (\cf \cite{Cisinski}) we also have $i$-th Chern class defined for elements of homotopy invariant $K$-theory satisfying that for any scheme $X$ the diagram
$$
\xymatrix{K(X)\ar[r]^-{c_i} \ar[d] &\EE^{2i,i}(X )\\
KH(X)\ar[ur]_-{c_i}}
$$ 
commutes.  We also obtain  higher Chern classes $c_{i,r}\colon KH_i(X)\to \EE^{2i-r,i}(X)$.
\end{parra}

\begin{parra}\label{Parra clases Chern Gysin KH y K}
Let $Z\to X$ be a closed immersion and denote $T=S^r\wedge X/X-Z$, then we denote
$$
c_{i,r}^Z\colon  KH_{Z,r} (X) := [T, \ZZ \times \Gr]\to [T, \Omega^\infty
\EE (i)[2i]] =: \EE^{2i-r,i}_Z(X).
$$
In addition, it follows from the construction that these Chern classes with support are functorial. As before, we obtain a commutative diagram
$$
\xymatrix{K_Z(X)\ar[r]^-{c_{i,r}^Z} \ar[d] &\EE_Z^{2i-r,i}(X )\\
KH_Z(X)\ar[ur]_-{c_{i,r}^Z}}
$$ 
for the $i$-th Chern class with support.  Let $i\colon Z \to X$ be a closed immersion, we denote by $i_\flat\colon \EE_Z(X)\to \EE(X)$ the natural morphism of forgetting support. We have  $i_\flat \circ (c_{i,r}^Z)= c_{i,r}$. We are ready to state the construction of the refined Gysin morphism obtained in \cite[\S 2]{Navarro advances}.

\end{parra}

\newpage

\begin{teo}
Let $(\EE, c_1)$ be an absolute oriented ring spectrum, there exist a unique way of assigning to any regular immersion a morphism $\bar i\colon \EE(Z)\to \EE_Z(X)$ such that:
\begin{enumerate}
\item
Normalization: if $i$ is of codimension one then $\bar i_* (a)= a\cdot c_1 ^Z (L_Z)$, where $L_Z$ denotes the line bundle given by the dual of the sheaf of ideals of regular functions vanishing on $Z$.

\item
Key formula: Consider the blow-up square
$$
\xymatrix{E\ar[r]^j\ar[d]_{\pi'}&B_ZX\ar[d]^\pi\\Z\ar[r]^i&X}
$$ 
and denote $K=\pi'^*N_{Z/X}/N_{E/B_ZX}$ the excess vector bundle. Then $ \pi^*(\bar i_*(a))=\bar j_*(\pi'^*(a))$ for all $a\in \EE(Z)$.
\end{enumerate}

\noindent
These Gysin morphism satisfies also the following properties:

\begin{enumerate}\setcounter{enumi}{2}
\item[2'.]
Excess intersection formula: Consider a cartesian square
$$
\xymatrix{ P \ar[r]^j \ar[d]_{\pi '} &X'\ar[d]^\pi \\
Z\ar[r]^i & X}
$$
where both $i$ and $j$ are regular immersions of codimension $n$ and $m$ respectively. Denote $K=\pi'^* N_{Z/X}/N_{P/X'}$ the excess vector bundle, then
$$
\pi^*\bar i_* (a)= \bar j_*(c_{n-m}(K)\cdot \pi '^*(a)) 
$$ 

\item
Functoriality: Let $j \colon Y\to Z$ be a regular immersion, then $(\overline{ij})_* =\bar i_* \bar j_*$.

\item
Projection formula: The Gysin morphism is $\EE(X)$-linear. In other words,
$$
a\cdot \bar i_* (b)= \bar i_* (i^* (a) \cdot b) \qquad \forall \ a \in \EE (X) \ , \ b \in \EE
(Z).
$$

\end{enumerate}
\end{teo}

The direct image $i_*$  in Thomason's algebraic $K$-theory satisfies the excess intersection formula as well as the normalization (\cf \cite{Thomason}). As a consequence, after a base change injective in cohomology the Gysin morphism is expressed in terms of Chern classes. Hence, applying an analogous argument as in \ref{Clases superiores con soporte} to $\KGL$ instead of $\EE$ and the normalization of the preceding result we deduce that for a regular immersion $i\colon Z\to X$ the square
$$
\xymatrix{K(Z)\ar[r]^{i_*}\ar[d] \ar[d] &K(X )\ar[d] \\
KH(Z)\ar[r]^-{ i_*}& KH(X)}
$$
also commutes. 

\begin{defi}\label{Defi Gysin}
Let $i\colon Z\to X$ be a regular closed immersion. We define the \textbf{Gysin morphism} to be $i_*\coloneqq i_\flat \circ \bar i \colon \EE(Z) \to \EE(X)$.  If $\MM$ is an absolute $\EE$-module, we define a refined Gysin morphism and its non-refined version in the $\MM$-cohomology by the formulas $\bar i_*(m)\coloneqq m \cdot \bar i_*(1)$ and  $i_*(m)\coloneqq i_\flat (\bar i_*(m))$, for all $m \in \MM(Z)$. 
\end{defi}

It is clear that the above theorem implies analogue formulas for the nonrefined Gysin morphism and also for modules. The reader may find a more elaborated proof of the following statement for modules and proper morphisms in \cite[2.5]{NN}.

\begin{coro}\label{Coro propiedades Gysin modulos}
Let $\EE$ be an absolute oriented ring spectrum and $\MM$ be an absolute $\EE$-module. The Gysin morphism, both in the $\EE$-cohomology and $\MM$-cohomology, satisfy the analogue of the properties 1, 2', 3 , 4 of the previous theorem.
\end{coro}
\qed

\begin{parra}\label{Propiedades clase de Thom}
In order to fix notation, we recall the theory of Thom classes. Let $V$ be a vector bundle of rank $n$ on a scheme $X$. We define the \textbf{Thom space} of $V$ as
$$
\Th(V)=V/V-0\simeq \bar V /\PP(V) \ \mbox{ in } \ \HH(X),
$$
where $\bar V=\PP(V\oplus 1)$. 
Its  $\EE$-cohomology fits into a long exact sequence
$$
\ldots \to \EE ^{**} (\Th (V)) \xrightarrow{\pi^*} \EE^{**}(\bar V)\to \EE^{**}(V)
\to \ldots
$$
where, from the projective bundle theorem, the third arrow is always a split epimorphism.  Let $0\to H\to \bar V \to \O_{\bar V}(1)\to 0$ be the dual of the canonical short exact sequence. We call the \textbf{Thom class} of $V$ to the class
$$
\t (V)\coloneqq c_n(H)=\sum_{i=0}^n (-1)^ic_i(V)\cdot x^{i} \in \EE^{2i,i}(\bar V),
$$
where $x=c_1(\O_{\bar V}(-1))$ and the equality follows from Cartan-Whitney formula. Since $\t (V)$ is zero in $\EE(\PP(V))$, we call the \textbf{refined Thom class} to the unique element
$$
\bar \t (V)\in \EE(\Th (V))\simeq \EE_X (\bar V)=\EE_X(V)
$$
such that $\pi^*(\bar \t (V) )=\t (V)$.

Any short exact sequence of vector bundles $0\to V'\to V\to V''\to 0$ induces an isomorphism 
$$
\EE(\Th(V'))\otimes \EE(\Th(V''))\xrightarrow \sim  \EE(\Th(V)).
$$
(\cf \cite[2.4.8]{Deglise}). This pairing satisfies that
\begin{equation}\label{Ecua mult clase Thom}
\bar \t(V')\otimes \bar\t (V'') \mapsto  \bar \t(V).
\end{equation}

Let $i\colon Z\to X$ be closed immersion between regular schemes. The deformation to the normal bundle gives an isomorphism $\EE_Z(X)\simeq \EE(\Th (N))=\EE_Z(N)$, where $N$ denotes the normal bundle and $Z$ is considered as a closed subscheme through the zero section. It follows from the unicity of direct image restricted to regular schemes that $ \bar i_*(1)=\bar \t(N)$. Let us now treat the general case:

\end{parra}

\begin{teo}\label{Teo clase de Thom}
Let $X$ be a finite dimensional noetherian scheme and  $V\to X$ be a vector bundle. Denote $ \bar V=\PP(V\oplus 1)$ the projective completion and $s\colon X\to \bar V$ the zero section. Then $\bar s_*(1)=\bar \t (V)$  in  $\EE(\Th(V))\simeq \EE_{X}(\bar V)$ or, equivalently,
$$
\bar s_*(1)= \bar \t (V) \ \mbox{ in } \   \EE(\bar V).
$$
\end{teo}
\demo  First observe that if $V$ is trivial the result follows from the normalization of the Gysin morphism. Now, applying Jouanolou's trick we can assume $X$ is affine. As a consequence $[V]=[V']$ in $K_0(X)$ if and only if they are stably equivalent. In other words, $V\oplus r\O\simeq V'\oplus r'\O$ for some $r, r'\in \NN$. Hence, any relation $[V]=[V']+ [V'']$  is given by a short exact sequence $0\to V'\oplus r'\O\to V\oplus r\O\to V''\oplus r''\O\to $ for some $r,r',r'' \in \NN$. Now, for any $r$ and any $V$ we have a natural short exact sequence $0 \to V \to V\oplus r\O \to r\O\to 0$. Thanks to the the additivity given by (\ref{Ecua mult clase Thom}) the statement $\bar \t (V+r\O)$, for any $r\in \NN$, is equivalent to that of $\bar \t (V)$. Hence, we can reduce to the case of line bundles.

Let $L$ be a line bundle, we have $\t (L)= c_1(L)-c_1(\O_{\bar L}(-1))$ and $s_*(1)=c_1(\I^*)$, where $\I$ stands for the sheaf of ideals of the zero section in $\bar L$.  This sheaf may be computed explicitly: the composition $\O_{\bar L}(-1) \to L\oplus \O \to L$ of the canonical morphism and the projection is an isomorphism out of the zero section, which induces  $ L^*\otimes \O_{\bar L}(-1)\simeq \I \to \O$.

Consider the canonical short exact sequence
$$
0\to \O_{\bar L}(-1) \to L\oplus \O \to Q\to 0
$$
where $Q$ is the canonical quotient bundle. Taking second exterior product it induces $L= \bigwedge^2 (L\oplus \O)=Q\otimes \O_{\bar L}(-1) $ so that $Q= L\otimes \O_{\bar L}(1)=\I^*$.

To finish the proof note that $c_1(L\oplus \O)=c_1(L)$. Also by the additivity of Chern classes $ c_1(L\oplus \O)=c_1(Q)+c_1(\O(-1))$ and we conclude.

\qed

\section{The two results}

Denote $P_q^d(\xi, c_1, \ldots c_{q-d}; c'_1, \ldots , c'_{q-d})$ the universal polynomial with integer coefficients defined in \cite[\S 1]{Jouanolou}. We are ready to prove the first result of this note.

\begin{teo}[Riemann-Roch without denominators]\label{Teo RR sin deno}
Let $i\colon Z\to X$ be a regular immersion of codimension $d$
. Then for any $q>0$ and any $a \in K_r(Z)$ we have in $ H_{\M,Z}^{2q-r}(X,\ZZ(q))$ that
\begin{equation}\label{Formula RR sin deno}
c_{q,r}^Z(\bar i(a))=\bar i(P_q^d(\mathrm{rk}(a), c_{1,r}(a), \ldots, c_{q-d,r}(a);
c_1(N_{Z/X}), \ldots ,  c_{q-d} (N_{Z/X}))).
\end{equation}
\end{teo}
\demo Thanks to the remark of \ref{Parra clases Chern Gysin KH y K} it is enough to prove the result for homotopy invariant $K$-theory. 

For convenience, denote $P_q^d(a, E)=P_q^d(\mathrm{rk}(a), c_{1,r}(a), \ldots, c_{q-d,r}(a); c_1(E), \ldots ,  c_{q-d} (E))$ for any vector bundle $E$ on $Z$. We consider the deformation to the projective closure of the normal bundle. That is to say, consider the commutative diagram
\begin{equation}\label{Deformacion a completado proyectivo}
\vcenter{\xymatrix{ \overline{N} \ar[r]^{i_0} & X' & X\ar[l]_{i_1}\\
Z \ar[r]\ar[u]^{s_0} & \AA ^1_Z \ar[u]^\iota & Z\ar[l] \ar[u]_i }}
\end{equation}
where $X'= B_{Y\times\{ 0\}}\AA ^1 _X$. For $U=X'-\AA^1_Z$, taking motivic cohomology in the deformation diagram we obtain
$$
\xymatrix{
& H_{\M}(U,\ZZ) & \\
H_{\M,Z}(\bar N,\ZZ) \ar[d]& H_{\M, \AA^1 _Z}(X',\ZZ)
\ar[l]_-{i_0^*} \ar[u]^h \ar[r]^-{i_1 ^*} & H_{\M,Z}(X,\ZZ) \\
H_{\M}(Z,\ZZ)\ar@<1ex>[u]^{\p_s}& H_{\M}(\AA^1_Z,\ZZ) \ar[l]_{\sim} \ar[u]^{\p_\iota}
\ar[r]^{v^*} &  H_{\M}(Z,\ZZ)\ar[u]_{ \p_i} }
$$
where $h=j^* \iota_\flat$ and $\bar s$ is injective (since it has a retract).

We now prove that if formula (\ref{Formula RR sin deno}) holds for $\iota\colon \AA^1_Z \to X'$ then it also holds for $i\colon Z \to X$. Since $v^*N_{\AA^1_Z/X'}=N_{Z/X}$ the refined versions of the excess intersection formula applied to the right square and the functoriality of higher Chern classes gives
$$
c_{q,r}^Z\bar i (v^*(b))=c_{q,r}^Z i_1^* (\bar \iota (b)) = i_1^* c_{q,r}^{\AA^1_Z} \q_\iota(b)
$$
and
$$
i_1 ^* (\bar {\iota}(P_q^d(b, N_{\AA^1_Z/X'})))= \bar {i}((P_q^d (v^*(b), N_{Z/X})))
$$
for $b\in KH(\AA^1_Z)$ such that $v^*(b)=a$. The last and the first term respectively are elements that theorem state that coincide.

Chasing the diagram we have that if $a \in H_{\M}(X',\ZZ)$ has $h(a)=0$ and $i_0^*(a)=0$ then $a=0 $. Since $h( c_{q,r}^{\AA^1_Z} \bar \iota (a)) = c_{q,r}j^* \bar \iota (a)=0 $ and $h \bar {\iota}=j^*\iota_*=0$ the last case left to prove is formula (\ref{Formula RR sin deno}) for the zero section $s\colon Z\to \bar N$ of the projective completion of the normal bundle. This is the case treated in the literature when $Z$ is smooth.

First recall that in \ref{Propiedades clase de Thom} we observed that $H_{\M}(\mathrm{Th} (N),\ZZ)\simeq H_{\M,Z}(\bar N,\ZZ)\to H_{\M}(\bar N,\ZZ)$ is injective, so it is enough to prove formula (\ref{Formula RR sin deno}) in $H_{\M}(\bar N,\ZZ)$.

We summarize in two lemmas computations which involve arguments in $K_0(\bar N)$ that do not need the smoothness assumption (\emph{cf.} \cite[4.3 and 4.4]{KY}).

\begin{lema}\label{Lema KY 1}
With the previous notations, denote $Q$ the canonical quotient bundle of $\bar N$ and $\t^{KH}(N)$ and $\t (N)$ the Thom class in $KH$ and $H_{\M}$ respectively. Then for any $b \in KH_r(\bar N)$ we have
$$
c_{q,r}(b \cdot \t^{KH}(N))=P_q^d(b, Q)\cdot \t (N).
$$
\end{lema}

\qed

\begin{lema}\label{Lema KY 2}
Let $p\colon \bar N\to Z$ be the projection. For any $a \in KH_r(Z)$ we have
$$
p^* (P_q^d(a, N))\cdot \t (N)=P_q^d(p^*a, Q)\cdot \t (N).
$$
\end{lema}

\qed

\medskip

To conclude recall that the fundamental class of the zero section coincide with the Thom class (Theorem \ref{Teo clase de Thom}). From here, the projection formula gives $s_* (a)=p^*( a)\cdot \t ^{KH}(N)$ and the analogous formula for motivic cohomology. With them, we conclude
$$
c_{q,r}(i_*(a))=c_{q,r}(p^*(a)\cdot \t (N))=i_*(P_q^d(a ,N)).
$$

\qed

\medskip

Recall that motivic cohomology is universal among cohomology theories such that Chern classes have an additive formal group law. Hence, we deduce the following result: 

\begin{coro}
Let $(\EE, c_1)$ be an absolute ring spectrum with an additive orientation. Let $i\colon Z\to X$ be a regular immersion of codimension $d$, then for any $q>0$ and any $a \in K_r(Z)$ we have in $\EE^{2q-r,q}(X)$ that
$$
c_{q,r}(i_*(a))= i_*(P_q^d(\mathrm{rk}(a), c_{1,r}(a), \ldots, c_{q-d,r}(a);
c_1(N_{Z/X}), \ldots ,  c_{q-d} (N_{Z/X}))).
$$
\end{coro}
\qed

\begin{parra}\label{Parra clases Chern relativas}
Recall that the $i$-th Chern class defines a map $c_i\colon \ZZ\times \Gr_X\to \Omega^\infty\EE (i)[2i]$ in $\HH_\bullet (X)$. Let $f\colon T\to X$ be a proper morphis, the $i$-th Chern class fits into a commutative diagram
$$
\xymatrix{\hofib_{\ZZ\times \Gr_X}(g)\ar[r]\ar@{-->}[d] &\ZZ\times
\Gr_X \ar[d]^-{c_i} \ar[r]&f_*f^*\ZZ\times \Gr_X \ar[d]\\
\hofib_{\EE}(g) \ar[r]& \Omega^\infty \EE (i)[2i]\ar[r]& f_*f^* \Omega^\infty \EE(i)[2i]}
$$
so that there exist a Chern class for the relative cohomology of $g$, which we still denote $c_i$. An analogous argument hold for the rank function. We deduce from the construction that these classes are functorial.

\end{parra}

One can check that in the universal polynomials $P_q^d(\xi, c_1, \ldots c_{q-d}; c'_1, \ldots , c'_{q-d})$ the elements $c'_1, \ldots , c'_{q-d}$ always appear multiplied by elements $\xi, c_1, \ldots c_{q-d}$ (\cf \cite[\S 3.2]{Levine} for example, where the reference denotes $Q_{d,q-d}$ instead of $P_q^d$). Hence, for any $m\in KH (f)$ and any vector bundle $m$ the element $P_q^d(m,E)$ belongs to $H_\M(f,\ZZ)$.

\begin{coro}\label{Coro RR sin deno relativo}
Let $i\colon Z\to X$ be a regular immersion of codimension $d$ and $g\colon T\to X$ be a proper morphism. Then for any $q>0$ and any $m \in K_r(f_{{\scriptscriptstyle Z}})$ we have in $ H_{\M,Z}(f,\ZZ)$ that
\begin{equation*}
c_{q,r}^Z(\bar i_*(m))=\bar i_*(P_q^d(\mathrm{rk}(m), c_{1,r}(m), \ldots, c_{q-d,r}(m);
c_1(N_{Z/X}), \ldots ,  c_{q-d} (N_{Z/X}))).
\end{equation*}
\end{coro}
\demo Taking fiber product of $T\to X$ along the deformation to the projective closure (\cf diagram (\ref{Deformacion a completado proyectivo})) we obtain the diagram
\begin{equation*}
\vcenter{\xymatrix{ \overline{N}_T \ar[r]^{i_0} & T' & T\ar[l]_{i_1}\\
T_Z \ar[r]\ar[u]^{s_0} & \AA ^1_{T_Z} \ar[u]^\iota & T_Z\ar[l] \ar[u]_i }}
\end{equation*}
where $T'=B_{T_Z\times\{0\}}\AA^1_T$ with proper morphisms $f'\colon T'\to X'$, $f_{{\scriptscriptstyle Z}}\colon T_Z\to Z$ and $\overline{f}\colon \overline{N}_T\to \overline{N}$. We obtain a commutative diagram
$$
\xymatrix{ & H_{\M}(g_{{\scriptscriptstyle
U}},\ZZ) & \\
H_{\M}(\bar g,\ZZ) \ar[d]& H_{\M}(g',\ZZ)
\ar[l]_-{i_0^*} \ar[u]^h \ar[r]^-{i_1 ^*} & H_{\M}(g,\ZZ) \\
H_{\M}(g_{{\scriptscriptstyle Z}},\ZZ)\ar@<1ex>[u]^{s_*}& H_{\M}(g_{\AA^1_Z},\ZZ)
\ar[l]_{\sim} \ar[u]^{\iota_*} \ar[r]^{v^*} & H_{\M}(g_{{\scriptscriptstyle
Z}},\ZZ).\ar[u]_{ i_*} }
$$
The arguments from Theorem \ref{Teo RR sin deno} transfer to the relative cohomology. Indeed, $H_\M(f_{\AA^1_Z},\ZZ)\simeq H_\M(f_{{\scriptscriptstyle Z}},\ZZ)$, the excess intersection formula holds for the Gysin morphism (Corollary \ref{Coro propiedades Gysin modulos}) Chern classes are functorial, and localization holds for relative cohomology. Hence, we only have to prove the formula for the direct image $s_*$.

Recall from Definition \ref{Defi Gysin} that the direct image in relative cohomology is defined by multiplying by the same fundamental class as in cohomology. In this case by the Thom classes by $\t^{KH}(N)$ and $\t(N)$. We conclude by recalling once again the computations from \cite[4.3 and 4.4]{KY}.

\qed

\begin{coro}
Let $(\EE, c_1)$ be an absolute ring spectrum with an additive orientation. Let $i\colon Z\to X$ be a regular immersion of codimension $d$ and $f\colon T\to X$ be a proper morphism. Then for any $q>0$ and any $m \in K_r(f_{{\scriptscriptstyle Z}})$ we have in $ \EE(f)$ that
\begin{equation*}
c_{q,r}(i_*(m))=i_*(P_q^d(\mathrm{rk}(m), c_{1,r}(m), \ldots, c_{q-d,r}(m);
c_1(N_{Z/X}), \ldots ,  c_{q-d} (N_{Z/X}))).
\end{equation*}
\end{coro}
\qed

\end{document}